\documentclass[10pt,a4paper, oneside]{article}
\usepackage{amsthm}
\usepackage[utf8]{inputenc}
\usepackage[T1]{fontenc}
\usepackage[centertags]{amsmath}
\usepackage{amsfonts}
\usepackage{graphicx}
\usepackage{amsmath}
\usepackage{amssymb}
\usepackage{latexsym}
\usepackage[mathscr]{euscript}
\usepackage{tikz,color,makeidx} 
\usetikzlibrary{calc,matrix,arrows,chains,positioning,scopes,decorations.pathmorphing,backgrounds,fit}

\numberwithin{equation}{section}
\newtheorem{teo}{Theorem}[section]
\newtheorem{cor}[teo]{Corollary}
\newtheorem{pro}[teo]{Proposition}
\newtheorem{lem}[teo]{Lemma}
\newtheorem{es}[teo]{Example}
\newtheorem{defi}[teo]{Definition}
\newtheorem{rem}[teo]{Remark}
\newtheorem{notat}[teo]{Notation}

\newcommand{\bdfn}{\begin{defi} \begin{rm}}
\newcommand{\edfn}{\end{rm} \end{defi}}
\newcommand{\bthm}{\begin{teo}}
\newcommand{\ethm}{\end{teo}}
\newcommand{\bprop}{\begin{pro}}
\newcommand{\eprop}{\end{pro}}
\newcommand{\bcor}{\begin{cor}}
\newcommand{\ecor}{\end{cor}}
\newcommand{\blem}{\begin{lem}}
\newcommand{\elem}{\end{lem}}
\newcommand{\bfact}{\begin{rem} \begin{rm}}
\newcommand{\efact}{\end{rm} \end{rem}}
\newcommand{\bex}{\begin{es} \begin{rm}}
\newcommand{\eex}{ \end{rm} \end{es}}
\newcommand{\bnot}{\begin{notat} \begin{rm}}
\newcommand{\enot}{\end{rm} \end{notat}}
\newcommand{\ten}{\otimes}
\newcommand{\quot}[2]{{\raisebox{.2em}{$#1$}\left/\raisebox{-.2em}{$#2$}\right.}}
\newcommand{\mv}{$\mathbf{ThMV} $}
\newcommand{\rmv}{$\mathbf{ThDMV^{'}} $}
\newcommand{\dashrmv}{$\mathbf{ThRMV^{'}}$}

\newcommand{\R}{\mathcal{R}}
\newcommand{\D}{\mathcal{D}}

\tikzset{node distance=2cm, auto}

\title{An analysis of the logic of Riesz Spaces with strong unit}

\author{ Antonio Di Nola,
Serafina Lapenta\\ 
{\small Department of Mathematics, University of Salerno,}\\ 
{\small  Via Giovanni Paolo II, 132 Fisciano (SA), Italy}\\
 {\small adinola@unisa.it, slapenta@unisa.it}
\and 
Ioana Leu\c stean \\
{\small Department of Computer Science,} \\
{\small Faculty of Mathematics and Computer Science, University of Bucharest,}\\
{\small Academiei nr.14, sector 1, C.P. 010014,  Bucharest, Romania}\\   
{\small ioana@fmi.unibuc.ro}
}
\date{}
\begin{document}

\maketitle

\begin{abstract}

We study \L ukasiewicz logic enriched with a scalar multiplication with scalars taken in $[0,1]$. Its algebraic models, 
called {\em Riesz MV-algebras}, are, up to isomorphism, unit intervals of Riesz spaces with a strong unit endowed with an appropriate structure. When only rational scalars are 
considered, one gets the class of {\em DMV-algebras} and a corresponding logical system. 
Our research follows two objectives. The first one is to deepen the connections between  functional analysis and the logic of Riesz MV-algebras. 
The second one is to study the finitely presented MV-algebras, DMV-algebras and Riesz MV-algebras, 
connecting them from   logical, algebraic and geometric perspective.  


\end{abstract}

\section*{Introduction}

In this paper we present the logical system  $\mathcal{R}$\L \ which extends the infinitely valued \L ukasiewicz logic with  a family of unary operators that are semantically interpreted as scalar multiplication with scalars from the real interval $[0,1]$. The category of the corresponding algebraic structures is equivalent with the category of Riesz spaces with strong unit. 

Recall that \L ukasiewicz logic \L\  is the system that has $\{\rightarrow, \neg\}$ as basic connectives and  whose axioms are L1-L4 below: 

(L1) $\varphi \rightarrow (\psi \rightarrow \varphi)$

(L2) $(\varphi \rightarrow \psi)\rightarrow((\psi \rightarrow \chi)\rightarrow(\varphi \rightarrow \chi))$

(L3) $ (\varphi \vee \psi)\rightarrow (\psi \vee \varphi)$

(L4) $(\neg \psi \rightarrow \neg \varphi)\rightarrow (\varphi \rightarrow \psi)$.

  The only deduction rule is {\em modus ponens}.

The corresponding algebraic structures, MV-algebras, were defined by C.C. Chang in 1958 \cite{Cha1}. Chang's definition was inspired by the theory of lattice ordered groups, consequently MV-algebras are structures $(A,\oplus, \neg, 0)$ satisfting some apropriate axioms, where   $x\to y=\neg x \oplus y$ for any $x,y$. The connection between MV-algebras and Abelian lattice ordered groups was fully investigated by D. Mundici \cite{Mun1} who proved the fundamental result that MV-algebras are categorically equivalent with Abelian lattice-ordered groups with strong unit. 

Since the standard model of \L\  is the real interval $[0,1]$ endowed with the \L ukasiewicz negation $\neg x=1-x$ and the \L ukasiewicz implication $x\rightarrow y=\min(1-x+y,1)$,  a natural problem was to study \L ukasiewicz logic enriched with a product operation, semantically interpreted in the real product on $[0,1]$. This line of research led to the definition of PMV-algebras, which are MV-algebras endowed with a internal binary operation, but in this case the standard model only generates a proper subvariety. Through an adaptation of Mundici's equivalences, PMV-algebras and their logic are connected with the theory of lattice-ordered rings with strong unit.

A different approach is presented in \cite{LeuMod,LeuRMV}, where the real  product on $[0,1]$ is interpreted as scalar multiplication, with scalars taken in $[0,1]$. The system $\mathcal{R}$\L\ further developed in this paper is a relatively simple extension of \L\ which is obtained by adding to the infinitely valued \L ukasiewicz logic a family of unary operators $\{\nabla_r\}_{r\in [0,1]}$, whose dual operators are semantically interpreted in a scalar multiplication. Consequently,  Riesz MV-algebras \--- the corresponding algebraic structures \--- are categorically equivalent with Riesz spaces with a strong unit. Note that our results are not the first connection between \L ukasiewicz logic and the theory of Riesz spaces, one can see \cite{BedeDiNola,EA} for previous investigations, while the seminal idea of a connection between Riesz Spaces and a subcategory of MV-algebras was given in \cite{DiNolaLettieri}.

In this paper, we establish connections between the system $\mathcal{R}$\L\ and elements of functional analysis, where Riesz spaces are fundamental structures. After some needed preliminaries, in Section \ref{seclogic} we define $\mathcal{R}$\L\ and prove some logic-related results, the main being a syntactical characterization of the uniform convergence in a Riesz Space. Using this concept of \emph{limit of formulas}, in Theorem \ref{imp} we describe any formula in $\mathcal{R}$\L\ as a sequence of formulas in the Rational \L ukasiewicz logic \cite{DMV} and in Section \ref{subsec:compl} we characterize two norm-completions of the Lindenbaum-Tarski algebra of $\mathcal{R}\L$. 

The central part of the paper deals with the problem of finding under which conditions we can establish categorical equivalences between subcategories of finitely presented MV-algebras, DMV-algebras and Riesz MV-algebras. This investigation is carried in Sections \ref{sec:tensor}-\ref{sec:polyandfp} using three different points of view, that take advantage of state-of-the-art techiques on MV-algebras: polyhedra, tensor product and the fresh approach via categories of presentations. Finally, we link all results together in Section \ref{L-theories}, where we use the syntactical notion of limit to fully describe those theories of $\mathcal{R}$\L \ that are axiomatized by formulas of \L.

\section{Preliminaries on algebraic structures}
{\em MV-algebras} are the algebraic counterpart of \L ukasiewicz propositional calculus. They are structures $(A,\oplus,^{*},0)$ of type $(2,1,0)$ such that, for any $x,y\in A$,
$$\begin{array}{ll}
 (A,\oplus,0) \mbox{ is an Abelian monoid}, & ({x^*})^* = x\\
 (x^* \oplus y)^* \oplus y = (y^* \oplus x)^* \oplus x, & 0^* \oplus x = 0^*.
	\end{array}$$
	
\noindent Further operations are defined as follows: $1$ is $ 0^*$, the  \L ukasiewicz
implication is $x\to y=x^* \oplus y$, the \L
ukasiewicz conjunction is $x\odot y =( x^*\oplus
y^*)^*$ and Chang's distance is $d(x,y)=(x^*\odot y)\oplus (x\odot y^*)$, for any $x$, $y\in A$.  If $x\vee y = x\oplus (y\odot x^*)$ and $x\wedge y = (x^*\vee y^*)^*$ then $(A,\vee,\wedge, 1,0)$ is a bounded distributive lattice.  If $x\in A$ and $n\in {\mathbb N}$ then $0x=0$ and $(n+1)x=(nx)\oplus x$. 

 The variety of MV-algebras is generated by  the standard structure
$$ [0,1]_{MV}=([0,1], \oplus, ^*,0),$$
 where $[0,1]$ is the real unit interval,  $x^*=1-x$ and  $x\oplus y=\min(1,x+y)$ for any $x,y\in [0,1]$. 
We urge the interested reader to consult  \cite{CDM,DinLeu}  for a basic introduction to MV-algebras and \cite{MunBook} for advanced topics.

A fruitful research direction has arisen from the idea of endowing MV-algebras with a product operation. In particular, when we consider a scalar multiplication, we obtain the notion of \textit{Riesz MV-algebras}.

A {\em Riesz MV-algebra} \cite{LeuRMV} is a structure 
$$(R, \oplus, ^*, 0, \{r\mid r \in [0,1]\})$$ such that $(R, \oplus, ^*, 0)$ is an MV-algebra and $\{r\mid r \in [0,1]\}$ is a family of unary operations such that  the following  properties hold for any $x,y \in A$ and $r,q\in [0,1]$:
$$\begin{array}{ll}
r (x\odot y^{*})=(r  x)\odot(r y)^{*},& 
(r\odot q^{*})\cdot x=(r x)\odot(qx)^{*},\\
r (q  x)=(rq) x,&
1 x=x.\end{array}$$

 The variety of Riesz MV-algebras is generated by 
 $$[0,1]_{\mathbb R}=([0,1],\oplus, ^*, \{\alpha\mid \alpha\in [0,1]\}, 0),$$
 where $([0,1],\oplus, ^*,0)$ is the standard MV-algebra and  $x\mapsto \alpha x$ is the real product of $\alpha$ and $x$, for $\alpha, x\in [0,1]$.

An intermediate class between MV-algebras and Riesz MV-algebras is the class of \textit{DMV-algebras} \cite{DMV}. By \cite[Corollary 2.1]{DMVNotes}, DMV-algebras can be seen as MV-algebras endowed with a scalar multiplication with scalars taken in $[0,1]_Q=[0,1]\cap \mathbb{Q}$. In this case the standard model is $[0,1]_{\mathbb Q}=([0,1]_Q, \oplus, ^*,\{q\mid q\in [0,1]_Q\})$, where $x\mapsto qx$ is the product of $q$ and $x$, for any $q\in [0,1]_Q$, and it generates the variety of DMV-algebras. Endowed with suitable morphisms, DMV-algebras and Riesz MV-algebras form categories that we shall denote by $\mathbf{DMV}$ and $\mathbf{RMV}$ respectively.

A fundamental result in the theory of MV-algebras is their categorical equivalence with Abelian lattice-ordered groups with a strong unit \cite{Mun1}. Such groups have compatible group and lattice structures, and a distinguished Archimedean element, called \textit{strong unit} \cite{Birk}. If $(G,u)$ is a lattice-ordered group with strong unit, then $([0,u],\oplus,^*, 0)$ is an MV-algebra, where $x\oplus y=(x+y)\wedge u$ and $ x^*=u-x$ for any $x,y\in [0,u]$. A similar correspondence is proved for Riesz MV-algebras and Riesz spaces with strong unit \cite{LeuRMV}, as well as for DMV-algebras and divisible $\ell u$-groups \cite{DMV}. Moreover, in \cite{DMVNotes}, the authors prove that DMV-algebras are categorical equivalent  with $\mathbb{Q}$-vector lattices with a strong unit, that is lattice-ordered linear spaces over $\mathbb{Q}$.

Finally, we mention that the structures considered in this paper are 
{\em semisimple}. Recall that semisimple MV-algebras are subdirect products of  MV-subalgebras of $[0,1]$ and they enjoy an convenient functional representation. Indeed, given a semisimple MV-algebra $A$, this isomorphic with a separating MV-subalgebra of $C(X)$, where $X$ is a suitable compact Hausdorff space \cite[Corollary 3.6.8]{CDM}. A DMV-algebra or a Riesz MV-algebra is semisimple if its MV-algebra reduct is semisimple.

\section{The logic $\mathcal{R}$\L\ of Riesz MV-algebras}\label{seclogic}
In this section we present the system $\mathcal{R}$\L\  \cite{LeuRMV}, a conservative extension of \L ukasiewicz logic \L,  that has  Riesz MV-algebras as models. 

Since Riesz MV-algebras are categorically equivalent with Riesz spaces with strong unit,  $\mathcal{R}$\L\ can be seen as the logical counterpart of the latter class of structures.

We review already known results, we add some new insight and we deepen the connection with notions from functional analysis. The system $\mathcal{R}$\L\ has uncountable syntax, but the compactness theorem still holds. We further define the {\em truth degree} and the {\em provability degree} of a formula in $\mathcal{R}$\L\ and we prove  Pavelka-style completeness. The system allows us to define a  syntactic notion of {\em limit}. Consequently, we express any formula of $\mathcal{R}$\L \ as a syntactic limit of formulas in $\mathcal{Q}$\L . Moreover,  in Section \ref{subsec:compl}  we show that the Lindenbaum-Tarski algebra of  $\mathcal{R}$\L  \  naturally becomes a metric space and we analyze its Cauchy completeness.

%We start with section with a remark.
%\bfact \label{rem:limit-rn}
%If $\lim_n r_n=r$ in $[0,1]$ then $\lim_n r_n x=r x$ in $A$. This easily follows by  $d(r_n x,rx)\leq d(r_n,r)x$ for any $r_n$, $r\in [0,1]$ and $x\in A$ and the fact that in any 
%arhimedean RMV-algebra we have $r_n\uparrow r$ implies $r_nx\uparrow rx$.
%\efact
We recall that the logical system $\mathcal R$\L, is obtained from \L ukasiewicz logic by adding a connective $\nabla_{r}$  for any $r\in [0,1]$. Hence the logical connectives are $\{\neg, \rightarrow\}\cup\{\nabla_r\mid r\in [0,1]\}$  and the axioms are the following:

(L) the axioms (L1)-(L4) of \L ukasiewicz logic, 

(R1) $\nabla_{r}(\varphi \rightarrow \psi)\leftrightarrow (\nabla_{r}\varphi \rightarrow \nabla_{r}\psi)$

(R2) $ \nabla_{(r \odot q ^*)}\varphi \leftrightarrow (\nabla_{q}\varphi \rightarrow \nabla_{r}\varphi)$

(R3) $\nabla_{r}(\nabla_{q}\varphi)\leftrightarrow \nabla_{r \cdot q} \varphi$

(R4) $\nabla_1 \varphi \leftrightarrow \varphi$,\\
The only deduction rule is \textit{Modus Ponens}. Note  that ``from $\varphi$, we prove $\nabla_r \varphi$ for any $r\in [0,1]$'' is a derivative rule. 

 The set of formulas of $\mathcal{R}$\L\  is  denoted $Form_{\mathcal{R}L}$.  For any $r\in [0,1]$ and $\varphi\in Form_{\mathcal{R}L}$ we set 
$$\Delta_{r}\varphi = \neg\nabla_{r}\neg\varphi.$$

If we only consider $\{\nabla_q \mid q\in [0,1]\cap \mathbb{Q}\}$, we obtain the system $\mathcal{Q}$\L \ that has DMV-algebras as models \cite{DMV, DMVNotes}. We refer to \cite{LeuRMV} and \cite{DMV,DMVNotes} for a more detailed account on the logics $\mathcal R$\L\ and $\mathcal{Q}$\L\ respectively.

\bfact [Completeness results]
\L ukasiewicz logic \L\  is complete with respect to its standard model $[0,1]_{MV}$ \cite{Cha}, the logic $\mathcal{R}$\L\ is complete with 
respect to $[0,1]_{\mathbb R}$ \cite{LeuRMV} and $\mathcal{Q}$\L\  is complete with respect to $[0,1]_{\mathbb Q}$ \cite{DMV}. 
\efact

Note that both $\mathcal{R}$\L\ and $\mathcal{Q}$\L\ are conservative extensions of \L .

As usual, a set of formulas $\Theta\subseteq Form_{\mathcal{R}L}$ is \textit{consistent} if there exists a formula $\varphi$ such that $\Theta \nvdash \varphi$. Moreover,  $\Theta$ is consistent if and only if  $e(\Theta)=\{1\}$ for some evaluation $e: Form_{\mathcal{R}L}\rightarrow [0,1]$. 

\bfact [Compacteness]
If $\Theta$ is a set of formulas such that for any finite $\Gamma\subseteq \Theta$ there exists a $[0,1]$-valuation 
$e_\Gamma$ such that $e_\Gamma(\Gamma)=\{1\}$, then there exists a $[0,1]$-evaluation $e$ such that $e(\Theta)=\{1\}$. The proof 
is a direct consequence of the fact that every finite subset of $\Theta$ is consistent.
\efact

For any formula $\varphi$ with $n$ variables we define the corresponding term function $f_\varphi:[0,1]^n\to [0,1]$ in the usual manner. 

Recall that a \textit{unital} piecewise linear function (PWL$_u$-function) with real (or integer, or rational) coefficients is a continuous function $f: [0,1]^n \rightarrow [0,1]$ such that there exist $f_1, \ldots f_m: \mathbb{R}^n \rightarrow \mathbb{R}$, $f_i= c_{i1}x_1 + \ldots + c_{in}x_n+b_i$ with $c_{ij}, b_i \in \mathbb{R}$ (or $c_{ij}, b_i \in \mathbb{Z}$, or $c_{ij}, b_i \in \mathbb{Q}$) and for any $\mathbf{x}=(x_1, \ldots , x_n)\in [0,1]^n$, there exists $i$ such that $f(\mathbf{x})=f_i(\mathbf{x})$.

One can easily see that $f_\varphi$ is a PWL$_u$-function with integer coefficients if $\varphi$ is a formula of \L , $f_\varphi$ is a PWL$_u$-function with rational coefficients if $\varphi$ is a formula of $\mathcal{Q}$\L \ and $f_\varphi$ is a PWL$_u$-function with real coefficients if $\varphi$ is a formula of $\mathcal{R}$\L.

 In \cite{McN,LeuRMV,DMV} the free MV-algebra, free DMV-algebra and free Riesz MV-algebra  are defined as algebras of $[0,1]$-valued functions. By \cite[Theorem 1 and Theorem 2]{McN}, $MV_n$ is the MV-algebra of PWL$_u$-functions in $n$ variables and integer coefficients; by \cite[Corollary 7]{LeuRMV}, $RMV_n$ is the Riesz MV-algebra of PWL$_u$-functions in $n$ variables and real coefficients; by \cite[Theorem 4.5]{DMV}, $DMV_n$ is the DMV-algebra of PWL$_u$-functions in $n$ variables and rational coefficients.

The Lindenbaum-Tarski algebras of \L, $\mathcal{Q}$\L\ and $\mathcal{R}$\L\ are denoted by $L$, $QL$ and $RL$, respectively. If $[\varphi]$ is the equivalence class of a formula $\varphi$  then, in the Lindenbaum-Tarski algebras $QL$ and $RL$, we define 
$r[\varphi]=[\Delta_r\varphi]$ for  any appropriate $r$.
 When only $n$ variables are considered, the Lindenbaum-Tarski algebras $L_n$, $QL_n$ and $RL_n$ are isomorphic to 
 $MV_n$, $DMV_n$ and $RMV_n$, respectively.

\subsection{Truth degree and provability degree}
\noindent In this section we define the {\em truth degree} and the {\em provability degree} of a formula, inspired by similar results in \cite{Hajek}. We also note that our system can be studied in the general framework developed in \cite{Novak}.

Let $\eta$ is a fixed theorem of $\mathcal{R}$\L . 
For any $r\in [0,1]$ we denote by $\eta_r$ the formula $\Delta_{r}\eta$. Thus, in the Lindenbaum-Tarski algebra of ${\mathcal R}$\L, $[\eta_r]=[\Delta_r \eta]=r[\eta]=r\mathbf{1}$. In the next proposition we will prove that this approach does not depend upon the choice of $\eta$.

\bprop \label{7.18}
The following properties hold for any $r$, $q\in [0,1]$:\\
(a) If $\vdash \tau$ and $\tau_r=\Delta_r\tau$, then $\vdash \tau_r \leftrightarrow \eta_r$,\\
(b) If $\vdash \tau$ and $\tau_r=\Delta_r\tau$, then $\Theta\vdash\Delta_{r}\eta\rightarrow \varphi$ if and only if $\Theta\vdash\Delta_{r}\tau \rightarrow \varphi$,\\
(c) $\vdash \neg\eta_{r}\leftrightarrow \eta_{r^{*}}$,\\
(d) $\vdash(\eta_{r}\rightarrow\eta_{q})\leftrightarrow \eta_{r\rightarrow q}$,\\
(e) $\vdash \Delta_{r}\eta_{q}\leftrightarrow \eta_{rq}$,\\
(f) $e(\Delta_{r}\varphi)=re(\varphi)$ where $\varphi$ is a formula and $e$ is an arbitrary evaluation,\\
(g) $e(\eta_{r})=r$, where $e$ is an arbitrary $[0,1]$-evaluation,\\
(h) $r\leq q$ iff $\vdash\eta_{r}\rightarrow\eta_{q}$.
\eprop
\begin{proof}
It is straightforward by completeness of $\mathcal{R}$\L .
\end{proof}
\bfact
When $r$ and $q$ are chosen in $[0,1] \cap \mathbb{Q}$, Proposition \ref{7.18}(c)-(d) are the bookkeeping axioms of Rational Pavelka logic \cite[3.3]{Hajek}. Proposition \ref{7.18}(e) can be also considered a bookkeeping axiom, since it shows that the logical constants act naturally  with respect to the scalar operation. Thus, we can think of $\mathcal R$\L \  as an axiomatic extension of a Pavelka propositional calculus.
\efact

\bdfn
Let $\varphi$ be an arbitrary formula of $\mathcal{R}$\L , we define:\\
(1) the {\em truth degree} of $\varphi$, by

$\parallel\varphi\parallel=\min\{e(\varphi)\mid e\text{ is a }[0,1]\text{-evaluation}\}$,\\
(2) the {\em provability degree} of $\varphi$, by

$\mid \varphi \mid=\max\{r \in [0,1]\mid \quad \vdash\eta_{r}\rightarrow \varphi\}$.
\edfn 

Note that, in general, provability degree and truth degree are defined by infimum and supremum. In this case they are indeed  minimum and maximum.

\bprop[Pavelka completeness]\label{7.21}
If $\varphi$ is an formula of $\mathcal{R}$\L , then $$\mid \varphi \mid =\parallel\varphi\parallel.$$
\eprop
\begin{proof}
We recall that if $\varphi$ is a formula in the variables $v_1, \dots, v_n$, there exists a continuous PWL$_u$-function with real coefficients such that, for any evaluation $e$, $e(\varphi)= f_{\varphi}(e(v_1), \dots , e(v_n))$. Since $f_{\varphi}$ is a continuous function over a compact set, Weierstrass's extreme values theorem ensures the existence of the minimum value. This entails the existence of a minimal truth-value $e_m(\varphi)$, and $\mid\varphi \mid $ is indeed a minimum. 

On the other direction, by completeness of $\mathcal R$\L , for any $r$ such that $\vdash \eta_r \rightarrow \varphi$, we have $r \le e(\varphi)$ for any evaluation $e$. In particular we have $r \le e_m (\varphi)$. Being $e_m$ the minimal evaluation, $e_m(\varphi)\le e(\varphi)$ for any other evaluation $e$, that is $\vdash \eta_{e_m(\varphi)}\rightarrow \varphi$ by completeness, and $\parallel \varphi \parallel$ is a maximum.
\end{proof}

%\noindent Finally,  let $\eta$ and $\tau$ be theorems of $\mathcal R$\L . Then, for any $\Theta\subseteq Form_{\mathcal{R}L}$ and for any $r\in [0,1]$, the following hold:\\
%(1) $\vdash\Delta_{r}\eta\leftrightarrow \Delta_{r}\tau$,\\
%(2) $\Theta\vdash\Delta_{r}\eta\rightarrow \varphi$ if and only if $\Theta\vdash\Delta_{r}\tau \rightarrow \varphi$.\\
%This means that the choice of the theorem $\eta$ in the beginning of our Pavelka approach does not affect the truth degree and the provability degree of a formula.

\subsection{A logical approach to limits}\label{sec:limit}
\noindent Since any formula $\varphi$ of $\mathcal{R}$\L\  corresponds to a continuous piecewise linear function with real coefficients $f_\varphi:[0,1]^n\to [0,1]$, we shall now express the uniform convergence of functions as a syntactic notion in $\mathcal{R}$\L.

We firstly prove the following.
\blem\label{ajutlimit}
Let $A$ be a compact topological space,
$f_n:A \to \mathbb{R}$, $n\in\mathbb{N}$,
be a sequence such that there exists a decreasing sequence $(g_n)_n$
of continuous functions such that
$\displaystyle{\bigwedge_n g_n=0}$ (where $\displaystyle{\bigwedge}$
denotes the pointwise infimum)
and $|f(x)-f_n(x)| \leq g_n(x)$ for every $n \in \mathbb{N}$ and
$x \in A$. Then $(f_n)_n$ uniformly converges to $f$.
\elem
\begin{proof}
Since the sequence $(g_n)_n$ is decreasing and all $g_n$'s are continuous, the sequence converges pointwise to the zero function. Being $A$ compact, by the Dini's theorem,
 the convergence is uniform. Thus, for every $\varepsilon >0$ there is $k \in
\mathbb{N}$ with $0 \leq g_n(x) \leq \varepsilon$ for every
$n \geq k$ and $x \in A$. By hypothesis, from this it follows that
$0 \leq |f(x)-f_n(x)| \leq g_n(x) \leq \varepsilon$ for
every $n \geq k$ and $x \in A$ and the claim is settled.
\end{proof}

Now we are ready to prove our main result.

\bthm\label{deflimit}
 Let $(\varphi_n)_n\subseteq Form_{\mathcal{R}L}$ and $\varphi\in Form_{\mathcal{R}L}$. Then the following are equivalent:\\
(i)  $(f_{\varphi_n})_n$ uniformly converges to $f_\varphi$,\\
(ii) for any $r\in [0,1)$ there exists $k$ such that $\vdash \eta_r \rightarrow (\varphi\leftrightarrow\varphi_n)$ for any $n\geq k$,\\
(iii) there is a decreasing sequence $([\psi_n])_n$ in $RL$ such that $\bigwedge_n [\psi_n]=0$ (pointwise) and  $d([\varphi_n], [\varphi])\le [\psi_n]$ for any $n$,\\
(iv)  there exists an increasing  sequence $(r_n)_n$ in $[0,1]$ such that $\bigvee_n r_n=1$ and $\vdash \eta_{r_n}\rightarrow (\varphi\leftrightarrow\varphi_n)$ for any $n$.
\ethm
\begin{proof} 
(i)$\Leftrightarrow$(ii) Using the correspondence between functions and formulas, the condition $\vdash \eta_r \rightarrow (\varphi\leftrightarrow\varphi_n)$ becomes 
$\mid f_{\varphi_n} - f_{\varphi}\mid \leq f_{\eta_r^*}$, which means
$\mid f_{\varphi_n} - f_{\varphi}\mid \leq 1-r$. The equivalence is obvious.\\
(i)$\Leftrightarrow$(iii) Recall that $\bigwedge_n [\psi_n]=0$ is the pointwise convergence to 0. In one direction the result follows from the fact that uniform convergence implies order convergence \--- see for example \cite[Theorem 16.2]{RS} \--- and the functions associated to $\eta_{r}$ are constant. In the other direction, we have that $\bigwedge_n f_{\psi_n}=0$  and 
 $|f_\varphi(x)-f_{\varphi_n}(x)|=d(f_\varphi, f_{\varphi_n})(x)\le f_{\psi_n}(x)$ for any $x$. Thus, the result follows by 
 Lemma  \ref{ajutlimit}.\\
(iv) $\Rightarrow$ (iii) Choosing $\psi_n = \eta_{r_n^*}$ the claim is settled.\\
(iii) $\Rightarrow$ (iv)   Let $s_n$ be $\sup\{ f_{\psi_n}(\mathbf{x})\mid \mathbf{x}\in \text{supp}(f_{\psi_n})\}$. By Dini's Theorem $\{ [\psi_n]\}_n$ uniformly converges to $\mathbf{0}$, and since the uniform convergence implies the order convergence, 
$ \bigwedge s_n = 0$ and $f_{d(\varphi_n,\varphi)} \leq f_{\psi_n}\leq s_n$. If we take $r_n=s_n^*$ then 
 $\bigvee r_n=1$ and $r_n\leq f_{(\varphi_n\leftrightarrow \varphi)}$, and the conclusion follows.
\end{proof}

\bfact [Uniform convergence and order convergence]
Condition (iii) from the above proposition is similar to the property of order-convergence\footnote{Cfr \cite[Theorem 16.1]{RS} for the analogous definition in Riesz Spaces.} in the Lindenbaum-Tarski algebra $RL$, but it is stronger. Indeed, when one deals with functions in Riesz Spaces, the infimum of the functions does not need to coincide with the pointwise infimum. Nonetheless, it is easily seen that if $\inf_n g_n=g$, then $g(x)\le \bigwedge_n g_n(x)$ for any $x$. Thus, $\bigwedge_n g_n(x)=0$ for any $x$ forces $g=\mathbf{0}$.  Similarly, $\bigvee_n g_n(x)=1$ for any $x$ forces $\sup_n g_n=\mathbf{1}$. In \cite[Exercise 18.14(i)]{RS} one can find an example of a sequence of functions in $C([0,1])$ for which the supremum function is identically 1, but the pointwise supremum is not.  Finally, we remark that in general the order-convergence and the uniform convergence do not coincide. An example for $C([0,1])$ can be obtained from \cite[Example 16.18]{RS}. Whence, our condition is, so to speak, the best possible to have a logical definition of uniform convergence.
\efact

\bdfn 
If $(\varphi_n)_n\subseteq Form_{\mathcal{R}L}$ and $\varphi\in Form_{\mathcal{R}L}$, then $\varphi$  is {\em the limit} of 
the sequence $(\varphi_n)_n$ , in symbols $\varphi=\lim_n\varphi_n$, if the equivalent conditions from Proposition \ref{deflimit} hold.
\edfn

\subsection{From Rational \L ukasiewicz logic to $\mathcal R$\L \ via limits}
\noindent We recall that a formula $\varphi$ of $\mathcal{R}$\L \ is a formula
of $\mathcal{Q}$\L \ if contains only symbols $\nabla_q$ with 
$q\in [0,1]\cap \mathbb{Q}$. In the following, we shall call such formulas {\em rational}.

A sequence $(\varphi_n)_n$ of formulas is {\em increasing} (resp. {\em decreasing}) if $\vdash \varphi_n\rightarrow\varphi_{n+1}$ (resp. $\vdash \varphi_{n-1}\rightarrow\varphi_{n}$). We shall write $g_n \uparrow g$ if $g$ is the limit of $g_n$'s and the sequence is increasing and $g_n \downarrow g$ if $g$ is the limit of $g_n$'s and the sequence is decreasing, where the $g_n$'s are  piecewise linear functions.

%\bnot
%For a monotone sequence $\{g_n\}_{n\in \mathbb{N}}$, we shall $g_n \uparrow g$ if $g$ is the limit of $g_n$'s and the sequence is increasing and $g_n \downarrow g$ if $g$ is the limit of $g_n$'s and the sequence is decreasing.
%\enot

\blem\label{lem:seqRL}
For any function $f$ in $RMV_m$ there exist an increasing sequence of functions $\{g_n\}_{n\in \mathbb{N}}$ and a decreasing sequence of functions $\{h_n\}_{n\in \mathbb{N}}$, both in $DMV_m$, such that $g_n \uparrow f$ and $h_n \downarrow f$.
\elem
\begin{proof}
Since $RMV_m$ is isomorphic with the algebra of term functions of the logic $\mathbb{R}\mathcal{L}$, we shall prove the result by structural induction. All the technical properties used for uniform convergence can be found in \cite{RS}.

If $f=\pi_i$, one of the projection functions, it is enough to take $g_n=h_n=\pi_i$ for any $n\in \mathbb{N}$.

If $f=l^*$, then there exists sequences $\{t_n\}_{n\in \mathbb{N}}$ and $\{s_n\}_{n\in \mathbb{N}}$ such that $t_n \uparrow l$ and $s_n \downarrow l$. Since $l^*=1-l$, it is easily seen that $1-t_n \downarrow 1-l$ and $1-s_n \uparrow 1-l$ and the claim is settled.
If $f=f_1\oplus f_2$, that there exist $\{t_n^1\}_{n\in \mathbb{N}}$ and $\{s_n^1\}_{n\in \mathbb{N}}$ such that $t_n^1 \uparrow f_1$ and $s_n^1 \downarrow f_1$ and $\{t_n^2\}_{n\in \mathbb{N}}$ and $\{s_n^2\}_{n\in \mathbb{N}}$ such that $t_n^2 \uparrow f_2$ and $s_n^2 \downarrow f_2$. Thus, $t_n^1\oplus t_n^2\uparrow f_1\oplus f_2$ and $s_n^1\oplus s_n^2\downarrow f_1\oplus f_2$, which settles the claim.

Finally, if $f=\Delta_r g$, there exist $\{h_n\}_{n\in \mathbb{N}}$ and $\{l_n\}_{n\in \mathbb{N}}$ such that $h_n \uparrow g$ and $l_n\downarrow g$. Moreover, it is possible to find an increasing sequence and a decreasing sequence of rational numbers $\{t_n\}_{n\in \mathbb{N}}$ and $\{s_n\}_{n\in \mathbb{N}}$ such that $t_n \uparrow r$ and $s_n \downarrow r$. Being $RMV_n$ a semisimple algebra, $\Delta_{s_n}h_n \uparrow \Delta_r g$ and $\Delta_{r_n}l_n \uparrow \Delta_r g$. By \cite{DMVNotes}, functions of the type $\Delta_q g$, with $q$ rational and $g\in DMV_m$, are term functions for $\mathbb{Q}$\L\ and the claim is settled. 
\end{proof}

\bthm\label{imp}
For any formula $\varphi$ of $\mathcal R$\L \ there exist an increasing sequence of rational formulas  $(\alpha_n)_n$ and a decreasing sequence rational formulas $(\beta_n)_n$  such that 
$\varphi=\lim_n\alpha_n=\lim_n\beta_n$.
\ethm
\begin{proof}
It is a consequence of Theorem \ref{deflimit} and Lemma \ref{lem:seqRL}.
\end{proof}

\subsection{Completions of the \mbox{Lindenbaum-Tarski} algebra} \label{subsec:compl}
\noindent  In the sequel we show that the Lindenbaum-Tarski algebra $RL_n$ becomes a normed space in a natural way and we analyse two different norm completions. We note that this problem can be also studied  for $MV_n$ and $DMV_n$, but  in the context of Riesz MV-algebras it is intimately connected with the theory of Banach lattices, as shown in the final remarks of this section.

\bdfn
 For any $[\varphi]$ in the  Lindenbaum-Tarski algebra $RL_n$ we define:\\
(un) $\|[\varphi]\|_u=\sup\{f_\varphi(\mathbf{x})|\mathbf{x}\in [0,1]^n\}$,\\
(in) $I([\varphi])=\int f_\varphi(\mathbf{x})d\mathbf{x}$,\\
where $f_\varphi: [0,1]^n \rightarrow [0,1]$ is the piecewise linear function associated to $\varphi$. 
\edfn

In any  Riesz MV-algebra $R$ it is possible to define the {\em unit seminorm} $\|\cdot\|_u:R\to [0,1]$  by $\|x\|_u=\inf\{r\in [0,1]\mid  x\leq r1\}$ for any $x\in R$ \cite{LeuRMV}.
One can easily see that $\|[\varphi]\|_u$  from (un) is the unit seminorm on $RL_n$ which, in this case, is actually a norm.

We recall that  a \textit{state} defined on an  MV-algebra $A$ is a function $s: A\to [0,1]$ such that $s(1)=1$ and $s(x\oplus y)=s(x)+s(y)$ for all $x,y\in A$ with  $x\odot y=0$. A state of $A$ is said to be \textit{faithful} if $s(x)=0$ implies $x=0$.
{\em States} were defined  in \cite{Mu} and they  are  generalizations of finitely additive probability measures on boolean algebras.  If $R$ is a Riesz MV-algebra, $s:R\to [0,1]$ is a state of $R$ if it is a state of its MV-algebra reduct. Moreover, $(R,s)$ is a normed space. If we define $\rho_s(x,y)=s((x\odot y^*)\oplus (x^* \odot y))$ for any $x$, $y\in R$, then $\rho_s$ is a pseudometric on $R$. We say that $(R,s)$ is {\em state-complete} if $(R,\rho_s)$ is a complete metric space.
From \cite[Theorem 3.4]{Mu}, we infer that  $I$ is a faithful state on $RL_n$ which attains rational values on rational formulas. 

Consequently,  $(RL_n,\|\cdot\|_u)$ and $(RL_n, I)$ are normed spaces.

\bprop Assume $(\varphi_n)_n\subseteq Form_{\mathcal{R}L}$ and $\varphi\in Form_{\mathcal{R}L}$. Then:\\
(a) $\lim_n\varphi_n=\varphi$ iff 
$\lim_n\|[d(\varphi_n,\varphi)]\|_u=0$,\\
(b) if $\lim_n\varphi_n=\varphi$ then $\lim_n I(\varphi_n)=I(\varphi)$\\
(c) $\vdash\neg\varphi$ iff  $I(f_{\varphi})=0$.
\eprop
\begin{proof}
(a) and (b) are straightforward from the fact that the domain of the corresponding term functions is $[0,1]^n$ and the convergence of $(f_{\varphi_n})_n$ to $f_{\varphi}$ is uniform.

(c) On one hand, we have that $\vdash \neg \varphi$ implies $e(\varphi)=0$ for any evaluation, which entails   $f_{\varphi}=\mathbf{0}$ and $I(f_{\varphi})=0$. On other hand, being $f_{\varphi}$ a non-negative and continuous function, $\int f_{\varphi}(\mathbf{x})d\mathbf{x}=0$ implies $f_{\varphi}=\mathbf{0}$. Since $e(\varphi)=f_{\varphi}(e(v_1), \dots, e(v_n))$, we get $e(\varphi)=0$ for any evaluation, which is equivalent to $\vdash \neg \varphi$ by completeness.
\end{proof}

In the following we characterize the Cauchy completions of the normed spaces $(RL_n,\|\cdot\|_u)$ and $(RL_n, I)$.

\bthm
 The norm-completion of the normed space $(RL_n, \|\cdot\|_u)$ is isometrically isomorphic with $(C([0,1]^n),\|\cdot\|_\infty)$.
\ethm
\begin{proof}
 The maximal ideal space $Max(RL_n)$ of $RL_n$ coincide with the maximal ideal space of its MV-algebra reduct therefore, by \cite[Theorem 4.16(iv)]{MunBook}, it is homeomorphic with $[0,1]^n$. We now use the fact that $RL_n$ is dense in $C([0,1]^n)$, as proven in \cite[Lemma 7.4]{EnzoReggio}. 
\end{proof}

\bcor\label{cor:approx}
For any continuous function $f:[0,1]^m\to [0,1]$ there exists a sequence of formulas $\{ \varphi_n\}_{n\in \mathbb{N}}\subseteq Form_{\mathcal{R}L}$ such that $f=\lim_n f_{\varphi_n}$.
\ecor

\bfact\label{mundici}
If $MV_n$ is the Lindenbaum-Tarski algebra of \L ukasiewicz logic, then $(MV_n, I)$ is a normed space. Let $\mu$ be the Lebesgue measure associated to $I$, assume $L^1(\mu)_u$ is the algebra of $[0,1]$-valued integrable functions on $[0,1]^n$ and  define $s_\mu(\hat{f})=I(f)$ for any $f\in L^1(\mu)_u$ \--- note that $\hat{f}$ is the class of $f$, provided  we identify two functions that are equal $\mu$-almost everywhere. Then
$(L^1(\mu)_u, s_\mu)$ is the Cauchy completion of $(MV_n,I)$ by   \cite[Theorem 16.7]{MunBook}, where the author mention that $L^1(\mu)_u$ is not just an MV-algebra, but it is equipped with a far richer structure than $(MV_n, I)$.   Moreover, being $\mu$ finite and being each function of $L^1(\mu)$ bounded and measurable, $L^1(\mu)_u$ coincide with the unit interval of $L^{\infty}(\mu)$, but the latter is endowed with a different norm.
 \efact

\bthm
The normed space $(L^1(\mu)_u, s_\mu)$ from Remark \ref{mundici} is a state-complete Riesz MV-algebra. Moreover,  $(L^1(\mu)_u, s_\mu)$ is, up to isomorphism, the norm-completion of the normed space $(RL_n, I)$ 
\ethm
\begin{proof}
 By \cite[Theorem 4 and Proposition 1]{LeuSCRMV}, $(L^1(\mu)_u, s_\mu)$ is a norm-complete Riesz MV-algebra.  We can safely identify $MV_n$ with an MV-subalgebra of $RL_n$ and we note that $(MV_n,I)$ is a normed subspace of $(RL_n,I)$.
 By \cite[Theorem 16.7]{MunBook},  the algebra $(L^1(\mu)_u, s_\mu)$ defined above is the Cauchy completion of $(MV_n,I)$. Since $f_\varphi$ is an integrable function for any formula $\varphi$ of $\mathcal R$\L, we infer that the Cauchy completion of $(MV_n,I)$ coincide with the Cauchy completion of $(RL_n,I)$. 
\end{proof}

Since Riesz spaces play a central r\^ole in functional analysis, the theory of Riesz MV-algebras has been connected with $C^*$-algebras, $M$-spaces and $L$-spaces. In particular,  the analogues of Kakutani's representation theorems are known for Riesz MV-algebras. 
  The interested reader is referred to \cite{LeuRMV, LeuSCRMV} for 
further details.
  Following these connections and using the functor $\Gamma_{\mathbb R}$, the norm-completion of $(RL_n, \|\cdot\|_u)$  is, up to isomorphism, the unit interval of an M-space. Similarly, the norm-completion of $(RL_n, I)$  is, up to isomorphism, the unit interval of an L-space.

\section{Dualities and adjunctions for finitely presented algebras}

\noindent We recall that an \textit{ideal} for an MV-algebra is a downward closed set, which is also closed with respect to the sum $\oplus$. Ideals in a Riesz MV-algebra (or a DMV-algebra) coincide with ideals of its MV-algebra reduct. A \textit{principal ideal} is an ideal generated by one element of the algebra. If $A$ is an MV-algebra and $a\in A$ then the principal ideal generated by $a$ is 
$$(a]=\{x\in A\mid x\leq na \mbox{ for some } n\in {\mathbb N}\}.$$

 An MV-algebra (a Riesz MV-algebra, a DMV-algebra) is \textit{finitely generated} if it is generated by a finite set of elements, and it is \textit{finitely presented} if it is the quotient of a free finitely generated MV-algebra (Riesz MV-algebra, DMV-algebra) by a principal ideal.  Note that free algebras and finitely presented algebras are semisimple, see \cite{CDM} for the case of MV-algebras, \cite{DMVNotes} for the case of DMV-algebras and  \cite{RMV-DiNola} for the case of Riesz MV-algebras.
  
Finitely presented structures are intimately connected with logic since, up to isomorphism, they are Lindenbaum-Tarski algebras corresponding to a finitely axiomatizable theories with finitely many variables \cite[Lemma 3.19 and Theorem 6.3]{MunBook}. 

In this section we focus on the algebraic approach, which provides more tools for our investigations towards finding equivalences, dualities or adjunction between subclasses of finitely presented MV-algebras, Riesz MV-algebras and DMV-algebras. We present three different approaches and, in the last section of the chapter, we interpret our results within the logical systems studied in the previous section.

\subsection{The tensor product and the finitely presented structures} \label{sec:tensor}

\noindent In the sequel we denote by $\mathbf{MV_{ss}}$,  $\mathbf{DMV_{ss}}$ and  $\mathbf{RMV_{ss}}$ the full subcategories of semisimple algebras. Our main tool is the tensor product of semisimple MV-algebras, defined by D. Mundici in \cite{Mun}, which allows us to complete the above diagram with adjoint functors. Let us now present few significant results.
 
Let $\ten$ be the tensor product of semisimple MV-algebras.  In \cite{LLTP1} is proved that $(\mathcal{U}_{\mathbb R}, \mathcal{T}_\ten)$ is an adjoint pair, where $\mathcal{U}_{\mathbb{R}}: \mathbf{RMV_{ss}} \rightarrow \mathbf{MV_{ss}}$ is the forgetful functor and $\mathcal{T}_{\ten}: \mathbf{MV_{ss}}\rightarrow \mathbf{RMV_{ss}}$ is defined on objects by $\mathcal{T}_{\ten}(A)=[0,1]\ten A$ and on morphisms through a suitable universal property: if $\iota_A:A\to \mathcal{T}_\ten (A)$ is the standard embedding, for any $A\in \mathbf{MV_{ss}}$, then for any $f\colon A\to B$, $\mathcal{T}_\ten (f)$ is defined as the unique map that satisfies $\iota_B\circ f=\mathcal{T}_\ten (f) \circ \iota_A$. An analogous result is proved in \cite{DMVNotes} for DMV-algebras but, in this case, the left adjoint functor $\mathcal{D}_{\ten}: \mathbf{MV_{ss}}\rightarrow \mathbf{DMV_{ss}}$ is defined by $\mathcal{D}_{\ten}(A)=([0,1]\cap \mathbb{Q})\ten A$. 
 
 With the above notations, define $\mathcal{T}_{\ten}^{\mathbb Q}: \mathbf{DMV_{ss}}\rightarrow \mathbf{RMV_{ss}}$  by 
 \begin{center}
 $\mathcal{T}_{\ten}^{\mathbb Q}(A)=\mathcal{T}_{\ten}({\mathcal U}(A))$ for any $A$ in $\mathbf{DMV_{ss}}$. 
\end{center}

We have the following diagram. 
 
\begin{center}
\begin{tikzpicture}

  \node (B) {$\mathbf{DMV_{ss}}$};
  \node (D) [right of=B] {$\mathbf{MV_{ss}}$};
   \node(G)[left of=B] {$\mathbf{RMV_{ss}}$};

   \draw[->] (B) to node [swap] {$\mathcal{U}$} (D);
   \draw[->] (G) to node [swap] {$\mathcal{U}$} (B);
 
 \draw[->,bend left] (D) to node {$\mathcal{T}_{\ten}$} (G);
 \draw[->,bend right] (D) to node [swap]{$\mathcal{D}_{\ten}$} (B);
 \draw[->,bend right] (B) to node [swap]{$\mathcal{T}_{\ten}^{\mathbb Q}$} (G);

\end{tikzpicture}
\end{center}

\bfact
The functors $\mathcal{T}_{\ten}^\mathbb{Q} \circ \mathcal{D}_\ten$  and $\mathcal{T}_\ten$ are naturally isomorphic, since they  are adjoints to the same forgetful functor.
\efact

In the following we  explore the behaviour of the tensor functors  $\mathcal{D}_\ten$ and $\mathcal{T}_\ten$ with respect to finitely presented algebras. 

We start with the simple result that an homomorphism of MV-algebras preserves principal ideals.   

\blem \label{lem:id1}
Let $I=(a]$ be a principal ideal in an MV-algebra $A$. For any homomorphism of  MV-algebras $f:A \to B$, $(f(I)]=(f(a)]$.
\elem
\begin{proof}
It is a routine application of the definitions of MV- homomorphisms and principal ideals.
%Let $x\in (f(I)]$. There exist $a_1, \ldots , a_m\in I$ and $n_1, \ldots n_m \in \mathbb{N}$ such that $x\le n_1 f(a_1)\oplus \ldots \oplus n_mf(a_m)= f(n_1a_1\oplus \ldots \oplus n_m a_m)$. Since $n_1a_1\oplus \ldots \oplus n_m a_m \in I$, there exists $k\in \mathbb{N}$ such that $n_1a_1\oplus \ldots \oplus n_m a_m \le ka$, then $x\le kf(a)$, i.e. $x\in (f(a)]$. Conversely, $f(a) \in f(I)$ implies $(f(a)]\subseteq (f(I)]$.
\end{proof}

 Let $P$ be either $[0,1]$ or $[0,1]_Q=[0,1]\cap \mathbb{Q}$ and $\iota_A: A\rightarrow P\ten A$ be the embedding in the semisimple tensor product. We recall that by  \cite[Proposition 5.1]{LLTP1}, $[0,1]\ten MV_n\simeq RMV_n$ and by \cite[Theorem 3.3]{DMVNotes}, $[0,1]_Q\ten MV_n \simeq DMV_n$.

\bthm\label{prop:finpres}
Let A be a finitely presented MV-algebra, and let $J$ be an  ideal of $MV_n$ such that $A\simeq MV_n/J$. Then
$$ P \ten (\quot{MV_n}{J})\simeq \quot{(P\ten MV_n)}{(\iota_A(J)]}.$$
\ethm
\begin{proof}
In the sequel, we will denote the embedding $\iota_A:A \hookrightarrow P\ten A$ with $\iota$.\\
We define the map $\beta_1$

$\beta_1:P \times \quot{MV_n}{J} \rightarrow \quot{(P\otimes MV_n)}{(\iota(J)]}$ by\\
$\beta_1(\alpha , [f]_J)=[\alpha \otimes f]_{(\iota(J)]}$. By definition of $\ten$, $\beta_1$ is trivially a bimorphism. Moreover, it is well defined: for $h\in [f]_J$ we prove that $[\alpha \ten h]_{(\iota(J)]}=[\alpha \ten f]_{(\iota(J)]}$. We have

$\alpha \ten h \in [\alpha \ten f]_{(\iota(J)]}$ iff $d(\alpha \ten h, \alpha \ten f)\in (\iota(J)]$ iff $d(\alpha \iota (h), \alpha \iota (f)))\in (\iota (J)]$.\\
Since $P\ten MV_n$ is either a Riesz MV-algebra or a DMV-algebra, by \cite[Lemma 3.19 and Proposition 3.10]{LeuMod} we get

$d(\alpha \iota (h), \alpha \iota (f)))=\alpha d(\iota(h),\iota(f))= \alpha (\iota (d(h,f)))$.\\
By definition $h\in [f]$ implies $d(h,f)\in J$, therefore $\iota (d(h,f))\in \iota(J) \subseteq (\iota(J)]$ and trivially $d(\alpha \ten h, \alpha \ten f)=\alpha (\iota (d(h,f)))\in (\iota(J)]$.

By universal property of the semisimple tensor product, there exists a homomorphism of MV-algebras $\lambda$

$\lambda : P\ten \quot{MV_n}{J}\rightarrow \quot{(P\otimes MV_n)}{(\iota(J)]}$ such that $\lambda (\alpha \ten [f]_J)=[\alpha \ten f]_{(\iota(J)]}$.\\
By \cite[Corollary 2]{LeuRMV} (or \cite[Lemma 3.1]{DMVNotes} in the case of DMV-algebras), $\lambda $ is a homomorphism of Riesz MV-algebras (or DMV-algebras).\\
On the other side, consider the bimorphism

$\beta_2: P\times MV_n\rightarrow P\ten \quot{MV_n}{J}$ defined by $\beta_2(\alpha, f)=\alpha \ten [f]_J$.\\
By universal property there exists a map $\Omega$

$\Omega: P\ten MV_n\rightarrow P\ten \quot{MV_n}{J}$ such that $\Omega (\alpha \ten f)=\alpha \ten [f]_J$.\\
By Lemma \ref{lem:id1}, $(\iota(J)]=(\iota(g)]=(1\ten g]$, therefore if $\mathbf{h}\in (1\ten g]$ there exists $n\in \mathbb{N}$ such that $\mathbf{h}\le n(1\ten g)$. Then $\Omega(\mathbf{h})\le n\Omega(1\ten g)=1\ten [g]_J$; since $g\in J$ the latter equals to $0$, therefore $\mathbf{h}\in Ker(\Omega)$ and $(\iota (J)] \subseteq Ker(\Omega)$. By the general homomorphism theorems \--- which hold for MV-algebras \--- there exists a map $\gamma$ 

$\gamma : \quot{P\ten MV_n}{(\iota(J)]}\rightarrow P\ten \quot{MV_n}{J}$  such that $\gamma ([\alpha \ten f]_{(\iota(J)]})=\Omega(\alpha \ten f)=\alpha \ten [f]_J$.\\
On generators we have
 
$(\gamma \circ \lambda)(\alpha \ten [f]_J)=\gamma ([\alpha \ten f]_{(\iota(J)]})=\alpha \ten [f]_J= \mathbf{I}(\alpha \ten [f]_J), $\\
and by universal property $\gamma \circ \lambda =\mathbf{I}$. Denoted by $\pi$ the canonical epimorphism $\pi: P\ten MV_n\rightarrow \quot{P\ten MV_n}{(\iota(J)]}$, we get

$ (\lambda \circ \gamma \circ \pi)(\alpha \ten f)= \lambda (\alpha \ten [f]_J)=[\alpha \ten f]_{(\iota(J)]}=(\mathbf{I}\circ \pi)(\alpha \ten f)$\\
and again by universal property, $\lambda \circ \gamma \circ \pi=\mathbf{I}\circ \pi$. Being $\pi$ is surjective, it follows that $\lambda \circ \gamma =\mathbf{I}$, and the claim is settled.
\end{proof}

\noindent If $A$ is a finitely presented MV-algebra, Theorem \ref{prop:finpres} ensures that $\mathcal{T}_{\ten}(A)$ is a finitely presented Riesz MV-algebra and $\mathcal{D}_{\ten}(A)$ is a finitely presented DMV-algebra.
We can therefore restrict the tensor functors to the full subcategories whose objects are finitely presented structures. 

\bfact \label{rem:idRMV}
Let $I$ be a principal ideal in $DMV_n$, generated by $a\in DMV_n\setminus MV_n$. By \cite[Remark 2.4]{LeuDia} and \cite[Theorem 33]{DMVNotes} there exists a natural number $k$ and $a_1,\ldots, a_k\in MV_n$ such that $a=\frac{1}{k}a_1\oplus\cdots \oplus\frac{1}{k}a_k$. Let $b=a_1\oplus\cdots \oplus a_k$. Since $\frac{1}{k}a_i \le a_i$, we have $a\le b$ therefore the ideal generated by $a$ is included in the ideal generated by $b$, in symbols $(a] \subseteq (b]$. For the other inclusion, $\delta_k(a_i)=\frac{1}{k}a_i\le a$ for any $i$, therefore $a_i=k\delta_k(a_i)\in (a]$. Hence $b\in (a]$,  and $(a] = (b]$. 

This means that we can replace the generator of a principal ideal in $DMV_n$ by an element of $MV_n$.
\efact

We denote $\mathcal{D}_{\ten}^{\mathbf{fp}}$ the restriction and co-restriction of  $\mathcal{D}_{\ten}$ to full subcategories of finitely presented structures.

\bthm \label{thm:algebrasD}
The functor $\mathcal{D}_{\ten}^{\mathbf{fp}}: \mathbf{MV_{fp}}\to \mathbf{DMV_{fp}}$ is faithful and essentially surjective.
\ethm
\begin{proof}
Let $D\in \mathbf{DMV_{fp}}$. There exists $n$ such that $A\simeq \quot{DMV_n}{I}$ where $I\subseteq DMV_n$ is a principal ideal. 
By Remark \ref{rem:idRMV} there exists $b\in MV_n$ such that $I=(b]$. By Theorem \ref{prop:finpres}, $\quot{DMV_n}{(b]_{DMV}} \simeq [0,1]_Q \ten \quot{MV_n}{(b]_{MV}}$, hence $A\simeq\mathcal{D}_{\ten}^{\mathbf{fp}}\left(\quot{MV_n}{(b]_{MV}}\right)$
and $\mathcal{D}_{\ten}^{\mathbf{fp}}$ is essentially surjective. 
Let $\sigma,\tau: A\to B$ arrows in $\mathbf{MV_{fp}}$ such that $\mathcal{D}_{\ten}^{\mathbf{fp}}(\sigma)=\mathcal{D}_{\ten}^{\mathbf{fp}}(\tau)$. In particular $\mathcal{D}_{\ten}^{\mathbf{fp}}(\sigma)(\iota_A(a))=\mathcal{D}_{\ten}^{\mathbf{fp}}(\tau)(\iota_A(a))$, and by definition of $\mathcal{D}_{\ten}^{\mathbf{fp}}$ on arrows, $\iota_B(\sigma(a))=\iota_B(\tau(a))$, for any $a\in A$. Being $\iota_B$ injective, we infer $\sigma=\tau$.
\end{proof}

To prove a similar result for finitely presented Riesz MV-algebras we need extra conditions of objects. Indeed, the next proposition shows that  we cannot prove the analogue of Remark \ref{rem:idRMV} in the case of Riesz MV-algebras.  We refer to Section \ref{sec:polyandfp} for the notion of a polyhedron and for suggestions for further readings.

\bprop\label{nonmv}
There exists at least a principal ideal in $RMV_n$ that is not generated by an element of $MV_n$.
\eprop
\begin{proof}
Let $P\subseteq [0,1]^n $ be a non-rational polyhedron. By \cite[Theorem 3.3]{RMV-DiNola} there exists $f\in RMV_n$ such that $P=f^{-1}(0)$. Consider now $I=(f]$, and assume that we can find $g\in I$ such that $g\in MV_n$ and $I=(g]$. Trivially, since $f\le kg$ and $g\le mf$ for appropriate $m,k\in \mathbb{N}$, $f(\mathbf{x})=0$ iff $g(\mathbf{x})=0$, for any $\mathbf{x}\in [0,1]^n$. This means $P=g^{-1}(0)$, which contradicts \cite[Proposition 2.1]{MarraMundici}, that is the fact that rational polyhedra are exactly zerosets of piecewise linear functions with integer coefficients (see also \cite[Lemma 2.5]{LucaEnzo}).
\end{proof}

A concrete example is the following. 
\bex \label{example:nonmv}
Consider the function $f\in RMV_1$ defined by $f(x)=-\sqrt{2}x+1$ for $x\le \frac{\sqrt{2}}{2}$ and $f(x)=0$ otherwise. Thus, $f^{-1}(0)=[\frac{\sqrt{2}}{2}, 1]$. Since $[\frac{\sqrt{2}}{2}, 1]$ has an irrational vertex, we can infer the desired conclusion as in Proposition \ref{nonmv}.
%If there exist a $g\in MV_1$ such that $(f]_{RMV}=(g]_{RMV}$, then in particular $f$ and $g$ have the same zerosets. Thus, $g^{-1}(0)=[\frac{\sqrt{2}}{2},0]$, which is an irrational polyhedron and this is a contradiction, since zerosets of elements of $MV_n$ are rational polyhedra by \cite[Lemma 3.2]{LucaEnzo}.
\eex

At this stage, we can only prove the following. 

\blem \label{thm:algebrasR}
The functor $\mathcal{T}_{\ten}^{\mathbf{fp}}: \mathbf{MV_{fp}}\to \mathbf{RMV_{fp}}$ is faithful.
\elem
\begin{proof}
It follows by the general properties of the tensor functor, as in the proof of Theorem \ref{thm:algebrasD}. 
\end{proof}

The previous result suggests that $\mathbf{MV_{fp}}$ is equivalent with a  subcategory of   $\mathbf{DMV_{fp}}$ that has the same objects as  $\mathbf{DMV_{fp}}$, but fewer morphisms. A similar remark holds if we consider a subcategory of $\mathbf{RMV_{fp}}$ instead of $\mathbf{DMV_{fp}}$. Thus,  we now define concrete subcategories of $\mathbf{DMV_{fp}}$ and $\mathbf{RMV_{fp}}$ that are equivalent to $\mathbf{MV_{fp}}$. 

\blem \label{morfcond}
Let $A$, $B$ be MV-algebras and $\iota_A:A\to {\mathcal D}_\ten (A)$, $\iota_B:B\to {\mathcal D}_\ten(B)$ the canonical embeddings of the tensor product. For 
a morphism $f_\ten: {\mathcal D}_\ten(A)\to {\mathcal D}_\ten(B)$ the following are equivalent: \\
(a) $f_\ten(\iota_A(A))\subseteq \iota_B(B)$,\\
(b) there exists a unique morphism $f:A\to B$ such that $f_\ten\circ\iota_A=\iota_B\circ f$. 
\elem
\begin{proof}
For the nontrivial implication, we set $f(a)=b$ if $a\in A$ and $b\in B$ such that $f_\ten(\iota_A(a))=\iota_B(b)$.  The uniqueness follows from the universal property of the tensor product. 
\end{proof}

Let $\mathbf{DMV_{fp}^{'}}$ be the  subcategory of $\mathbf{DMV_{fp}}$  whose objects are algebras $D={\mathcal D}_\ten(A)$  where $A$ is a finitely presented MV-algebra and whose morphisms are  $g:
{\mathcal D}_\ten(A)\to {\mathcal D}_\ten(B)$ such that    
$g(\iota_A(A))\subseteq \iota_B(B)$.  Similarly, let $\mathbf{RMV_{fp}^{'}}$ be the  subcategory of $\mathbf{RMV_{fp}}$  whose objects are algebras $R={\mathcal T}_\ten(A)$  where $A$ is a finitely presented MV-algebra and whose morphisms are  $g:
{\mathcal T}_\ten(A)\to {\mathcal T}_\ten(B)$ such that    
$g(\iota_A(A))\subseteq \iota_B(B)$.

\bprop\label{preq}
The categories $\mathbf{MV_{fp}}$, $\mathbf{DMV_{fp}^{'}}$ and $\mathbf{RMV_{fp}^{'}}$ are equivalent.
\eprop 
\begin{proof}
One can easily see that the co-restrictions of ${\mathcal D}_\ten$ and 
${\mathcal T}_\ten$ are full and  essentially surjective. The faithfulness follows by Lemma \ref{morfcond}, which is easily adapted to the case of Riesz MV-algebras. 
\end{proof}

Note that, as a consequence of Proposition \ref{nonmv}, the functor $\mathcal{T}^{fp}_{\ten}:\mathbf{MV_{fp}}\to \mathbf{RMV_{fp}}$ is not essentially surjective, while $\mathbf{RMV_{fp}^{'}}$ has exactly those objects that make the functor essentially surjective.

Let $\mathbf{DMV_{fp}^c}$ be the full subcategory of $\mathbf{DMV}$  whose objects are algebras $D={\mathcal D}_\ten(A)$  with $A$ being a finitely presented MV-algebra. Hence $\mathbf{DMV_{fp}^c}$ is equivalent to $\mathbf{DMV_{fp}}$ and 
$\mathbf{DMV_{fp}^{'}}$ is a subcategory equivalent to $\mathbf{MV_{fp}}$.  Note that $\mathbf{DMV_{fp}^{'}}$ and $\mathbf{DMV_{fp}^{c}}$ have the same objects, but $\mathbf{DMV_{fp}^{'}}$ has fewer morphisms.  
We say that the category $\mathbf{DMV_{fp}^{c}}$ is {\em concrete} since the objects are defined fixing a particular representation.  An equivalent concrete category can be  defined, by considering the full subcategory of $\mathbf{DMV}$ whose  objects have the  form $DMV_n/I$, where $I$ is a principal ideal \cite{LeoLuca}. Another concrete category  is defined in Section  \ref{sec:polyandfp} using polyhedra. 

Note that, in our setting, $A$ and $\iota_A(A)$ are isomorphic MV-algebras. Since $\iota_A(A)$ is an MV-subalgebra of an enriched structure - DMV-algebra or Riesz MV-algebra - we say that the elements of $\iota_A(A)$ are {\em MV-elements}. In order to get a category that is equivalent to $\mathbf{MV_{fp}}$ we have to select fewer morphisms, exactly those morphisms that  preserves the MV-elements.

\subsection{Categories of presentations}\label{pres}  
\noindent We proved that the category $\mathbf{MV_{fp}}$ of finitely presented algebras is equivalent to subcategories of DMV-algebras and Riesz MV-algebras whose objects are defined choosing a concrete representation. Fixing a representation might not be desirable when analysing abstract properties  therefore, in this section,  we present a different construction  inspired by the categorical approach to presentations and theories in the general framework of institutions \cite{Dia}. We analyse  finitely presented structures from the point of view of their presentations. In other words, we define categories whose objects are theories determined by finite presentations. 

We denote by $MV_X$, $DMV_X$ and $RMV_X$ the Lindenbaum-Tarski algebras over the set $X$ of variables, corresponding to the systems \L, $\mathcal Q$\L \ and $\mathcal R$\L \ from Section \ref{seclogic}. 
Let $\iota_X$ be the unique morphism of MV-algebras $\iota_X : MV_X \rightarrow DMV_X$ defined by $\iota_X(x)=x$ for any $x\in X$. $\iota_X$ is an trivially an embedding. In the following, we will often identify $S\subseteq MV_X$ with  $\iota_X(S)$. In particular, $MV_X$ will be identified with an MV-subalgebra of $DMV_X$.  Similarly, we shall identify $MV_X$ with an MV-subalgebra of $RMV_X$. In this approach, the elements of the isomorphic copy of $MV_X$ are the MV-elements of $DMV_X$ or $RMV_X$. 

The  theories of  \L \ are in  bijective correspondence with  the filters of $MV_X$ \cite[Proposition 6.3.15]{DinLeu} and, consequently,  with the ideals of $MV_X$. Thus, we represent a theory by a pair $(X,I)$, where $I$ is an ideal of $MV_X$. In the following we consider only finitely generated theories, but the development can be generalized to arbitrary ones.

\bdfn
We define the category \mv \ in the following way:
the objects are pairs $(X,I)$ with X a set and $I$ an finitely generated ideal of $MV_X$; 
the arrows are maps $\sigma:(X,I)\rightarrow (Y,J)$ such that $\sigma: MV_X\rightarrow MV_Y$ and  $\sigma (I) \subseteq J$.
\edfn

\bfact (1) In the general theory of institutions, morphisms of theories are usually chosen to be maps from variables to variables, expressing the fact that the truth is invariant under renaming. We use substitutions as morphisms, instead of using renaming maps.\\
(2) If $(X,I)$ and $(Y,J)$ from \mv \ are isomorphic then $\quot{MV_X}{I}$ and $\quot{MV_Y}{J}$ are isomorphic MV-algebras. Moreover, if $\quot{MV_X}{I}$ and $\quot{MV_Y}{J}$ are isomorphic MV-algebras then there are  morphisms $\sigma:(X,I)\to (Y,J)$ and  $\tau:(Y,J)\to (X,I)$ in \mv such that $\sigma(\tau(y))/J=y/J$ and $\tau(\sigma(x))/I=x/I$ for any $x\in X$ and $y\in Y$. \\
(3) In \cite{LeoLuca}, a different category of presentations is defined, by considering morphisms between the quotients $\quot{MV_X}{I}$ and $\quot{MV_Y}{J}$. This category is equivalent with $\mathbf{MV_{fp}}$ and they are obviously equivalent with a quotient category of \mv. 
%One can see  \cite{CMS} for more details on this subject. 
\efact

We define $\mathbf{ThDMV}$ in a similar way and we denote \rmv \  the subcategory of $\mathbf{ThDMV}$ that has the same objects but fewer morphims:  an arrow $\sigma:(X,I)\rightarrow (Y,J)$ from $\mathbf{ThDMV}$  is an arrow in \rmv \ if $\sigma (X)\subseteq MV_Y$. Note that we are asking for arrows that preserve MV-elements.
 
The following definition is essential for our development.

\bdfn\label{def:mvgen}
We say that an ideal $J$ of $DMV_X$ is \textit{MV-generated} if $J=(J\cap MV_X]_{DMV}$. We remark that the inclusion $(J\cap MV_X]_{DMV}\subseteq J$ is always satisfied. The MV-generated ideals of $RMV_X$ are similarly defined. Thus, an ideal is an MV-ideal if it is generated by its MV-elements and it is easily seen that a principal ideal of either $DMV_X$ or $RMV_X$ is MV-generated if and only if we can replace its generator with an element of $MV_X$.
\edfn

\blem \label{lem:idealsRMV}
Let $I$ be a principal ideal of $DMV_X$. Then $I$ is MV-generated and $\iota_X^{-1}(I \cap MV_X)$ is  a principal ideal in $MV_X$.
\elem
\begin{proof}
Let $I=(g]_{DMV}$. By Remark \ref{rem:idRMV}, there exists $g^*\in MV_X$ such that $I=(g^*]_{DMV}$. Hence $g^*\in I\cap MV_X$, $I\subseteq (I\cap MV_X]_{DMV}$ and $I$ is MV-generated.\\
Since the proof that $\iota_X^{-1}(I \cap MV_X)$ is an ideal of $MV_X$ is routine, let us prove that it is a principal ideal. By Remark \ref{rem:idRMV}, it is enough to consider the case $g\in I\cap MV_X$. We have $\iota_X^{-1}(g)\in \iota_X^{-1}( I\cap MV_X)$, therefore $(\iota_X^{-1}(g)]_{MV}\subseteq \iota_X^{-1}( I\cap MV_X)$. On the other side, for any $a \in \iota_X^{-1}(I\cap MV_X)$, we have $\iota_X(a) \in I\cap MV_X \subseteq I$, therefore there exists $n\in \mathbb{N}$ such that $\iota_X(a) \le ng$. This entails that $a\le n \iota_{X}^{-1}(g)$ and $\iota_X^{-1}(I\cap MV_X) \subseteq (\iota_{X}^{-1}(g)]_{MV}$ and the claim is settled.
\end{proof}

 \bthm \label{teo:catEq}
We define a functor $\mathcal{DT}:$ \mv $\rightarrow$ \rmv \ as follows: 
\begin{itemize}
\item \textit{Objects}: $\mathcal{DT}(X,I)=(X, (I]_{DMV})$, where 
$(I]_{DMV}$ is the $DMV$-ideal generated by $I$ in $DMV_X$; 
\item \textit{Arrows}: if $\sigma :(X,I)\rightarrow (Y,J)$ is an arrow in \mv , $\mathcal{DT}(\sigma)$ is the unique morphism $\sigma_D : (X, (I]_{DMV})\rightarrow (Y , (J]_{DMV})$ such that $\sigma_D: DMV_X \rightarrow DMV_Y$ and $\sigma_D (x)=\sigma(x)$ for any $x\in X$. 

\end{itemize}

The functor $\mathcal{DT}$ establishes a categorical equivalence between \mv \ and \rmv .
\ethm
\begin{proof}
We will prove that $\mathcal{DT}$ is full, faithful and essentially surjective.

\textit{$\mathcal{DT}$ is faithful}: for any $\tau, \sigma:(X, I)\rightarrow (Y, J)$, let $\tau_D = \sigma_D$. Then we have $\iota_Y(\tau(x))=\tau_D(\iota_X(x))=\sigma_D(\iota_X(x))=\iota_Y(\sigma(x))$. Being $\iota_Y$ injective, we get $\tau(x)=\sigma(x)$ for any $x\in X$. Since $X$ is the generating set, $\tau=\sigma$.

\textit{$\mathcal{DT}$ is full}: let $\tau: (X, (\iota_X(I)]_{DMV})\rightarrow (Y, 
(\iota_Y(J)]_{DMV})$ a morphism in \rmv . Then $\tau : DMV_X \rightarrow DMV_Y$, $\tau (X) \subseteq MV_Y$ and $\tau((\iota_X(I)]_{DMV}) \subseteq (\iota_Y(J)]_{DMV} $.\\
We define $\sigma : (X,I)\rightarrow (Y,J)$ as the unique map $\sigma: MV_X \rightarrow MV_Y$ such that $\sigma(x)=\tau(x)\in MV_Y$ for any $x\in X$. More precisely, $\iota_Y(\sigma(x)) = \tau (\iota_X(x))$, and since $X$ is the generating set, $\iota_Y \circ \sigma=\tau \circ \iota_X$. 
\begin{center}
\begin{tikzpicture}
  \node (A) {$MV_X$};
  \node (B) [right=1.5cm of A] {$MV_Y$};
  \node (C) [below of=A] {$DMV_X$};
  \node (D) [below of=B] {$DMV_Y$};
  \draw[->] (A) to node {$\sigma $} (B);
  \draw[->] (A) to node [swap]{$\iota_X$} (C);
  \draw[->] (C) to node [swap] {$\tau$} (D);
  \draw[->] (B) to node {$\iota_Y$} (D);
 \end{tikzpicture}
\end{center}
We have to prove that $\sigma (I) \subseteq J$. By hypothesis we get $\iota_Y(\sigma(I))= \tau (\iota_X(I))\subseteq (\iota_Y(J)]_{DMV}$.
Let $a \in I$, then $\sigma(a) \in \sigma(I)$ and $\iota_Y(\sigma(a))\in (\iota_Y(J)]_{DMV}$. Whence there exist $\alpha_1, \ldots , \alpha_m \in \mathbb N$ and $b_1, \ldots , b_m \in J$ such that $\iota_Y(\sigma(a))\le \delta_{\alpha_1}( \iota_Y(b_1))\oplus \ldots \oplus \delta_{\alpha_m}( \iota_Y(b_m)) \le \iota_Y(b_1)\oplus \ldots \oplus \iota_Y(b_m)$, since $\delta_k(x)\le x$ for any $k\in \mathbb N$.
Then $\iota_Y(\sigma (a))\le \iota_Y(b_1 \oplus \ldots \oplus b_m)$ implies $\sigma(a) \le b_1 \oplus \ldots \oplus b_m \in J $ so $\sigma(I) \subseteq J$ and $\tau=\mathcal{DT}(\sigma)$.

\textit{$\mathcal{DT}$ is essentially surjective}: we have to prove that for any $(Y, J)$ in \rmv \ there exists $(X,I)$ in \mv \ such that $\mathcal{DT}(X,I)=(Y,J)$. If we choose $I=\iota_X^{-1}(J \cap MV_X)$, the conclusion follows by Lemma \ref{lem:idealsRMV}. 
\end{proof}

The same approach works for theories in $\mathcal R$\L, \ but we have to restrict our attention to MV-generated ideals.  Indeed, by Proposition 
\ref{nonmv}, there are ideals in $RMV_X$ that are not MV-generated. 
Let \dashrmv \ be  the category whose objects are couples $(X,I)$, where $I\subseteq RMV_X$ is an ideal generated by $f\in DMV_X$ and whose morphisms $\tau : (X,I) \rightarrow (Y, J)$ are maps $\tau:RMV_X \rightarrow RMV_Y$ such that $\tau(I)\subseteq J$ and $\tau(X)\subseteq MV_Y$. 

\bprop \label{pro:eqPres}
The categories \mv , \rmv , and \dashrmv \  are equivalent.
\eprop
\begin{proof}
It is similar with the proof of  Theorem \ref{teo:catEq}. 
\end{proof}

Finally, we give the last point of view on the relations between finitely presented objects. 

\subsection{Polyhedra and finitely presented algebras}\label{sec:polyandfp}
In order to give out last characterization we need some preliminary notions of polyhedral geometry. An $m$-simplex in $\mathbb{R}^n$ is the convex hull $C$ of $m+1$ affinely independent points $\{ v_o, \ldots , v_m\}$ in the euclidean space $\mathbb{R}^n$; the points $v_i$ are called vertexes and $C$ is \emph{rational} if any vertex has rational coordinates.  A (\textit{rational}) \textit{polyhedron} is the union of finitely many (rational) simplexes.

If $P \subseteq [0,1]^n$ and $Q\subseteq [0,1]^m$ are polyhedra, then a map $z: P\rightarrow Q$ is a {\em $\mathbb{Z}$-map} if  $z=(z_1, \ldots z_m)$ where all $z_i:[0,1]^n\rightarrow [0,1]$ are piecewise linear functions with integer coefficients. Note that {\em $\mathbb Q$-maps} and {\em $\mathbb R$-maps} are  similarly defined, each corresponding to the given ring of coefficients.  

Let us denote by $\mathbf{Pol}_{[0,1]}^{\mathbb R}$ is the category of  polyhedra with $\mathbb{R}$-maps as morphisms, by $\mathbf{RatPol}_{[0,1]}^{\mathbb Q}$ the category of rational polyhedra with $\mathbb Q$-maps as morphisms and by $\mathbf{RatPol}_{[0,1]}^{\mathbb Z}$ the category of rational polyhedra with $\mathbb Z$-maps as morphisms, where all polyhedra are lying in a unit cube. Following the ideas of Baker and Beynon \cite{Baker,Bey,Bey1}, in  \cite{LucaEnzo,CDM} the authors proves that the categories $\mathbf{RatPol}_{[0,1]}^{\mathbb Z}$ and $\mathbf{MV_{fp}}$ are dual. This result is generalized to DMV-algebras in \cite{DMVNotes}, meaning that the categories $\mathbf{RatPol}_{[0,1]}^{\mathbb Q}$ and $\mathbf{DMV_{fp}}$ are dual, while in \cite{RMV-DiNola} it is proved that $\mathbf{Pol}_{[0,1]}^{\mathbb R}$ and $\mathbf{RMV_{fp}}$ are dual. 

We shall now further analyze these connections. Firstly, let us recall the definitions of the functors that define these dualities. 

Let $A$ be a finitely presented Riesz MV-algebra and assume $A\simeq \quot{RMV_n}{I}$, where $I$ is a principal ideal generated by $f\in RMV_n$. If $P=f^{-1}(0)$, then $P$ is a polyhedron in  $\mathbf{Pol}_{[0,1]}^{\mathbb R}$ and $A\simeq RMV_n\mid_P$, where $RMV_n\mid_P=\{g\mid_P\,\, \mid \,g\in RMV_n\}$.  Thus, in the case of polyhedra with $\mathbb{R}$-maps and finitely presented Riesz MV-algebras, the functor $\R : \mathbf{Pol}_{[0,1]}^{\mathbb R}\to\mathbf{RMV_{fp}}$ is defined as follows:

$\R (P)=RMV_n\mid_{P}$, for any  polyhedron $P\subseteq [0,1]^n$;

$\R (\lambda): \R (Q) \rightarrow \R (P)$ defined by $\R (\lambda)(f)=f \circ \lambda$, for any $\lambda: P \rightarrow Q$, with $P\subseteq [0,1]^n$ and $Q\subseteq [0,1]^m$.

The functors ${\mathcal M} : \mathbf{RatPol}_{[0,1]}^{\mathbb Z}\to\mathbf{MV_{fp}}$  and $\D : \mathbf{RatPol}_{[0,1]}^{\mathbb Q}\to\mathbf{DMV_{fp}}$ are defined similarly.  

This is the second concrete category of finitely presented structures. The first concrete representation for was given at the end of Section \ref{sec:tensor}.

\blem \label{teo:subDMV}
Let $\sigma:DMV_n\mid_P\to DMV_m\mid_Q$ be a map between finitely presented DMV-algebras, where  $P,Q$ are rational polyhedra  and assume $\lambda=(\lambda_1, \dots, \lambda_n):Q \to P$ such that ${\mathcal D}(\lambda)=\sigma$. Then $\sigma(\{\pi_1, \dots \pi_n \})\subseteq MV_m$ iff $\lambda_i \in MV_m$.
\elem
\begin{proof}
The result follows from the fact that $\sigma = \D (\lambda)$ implies   $\sigma(\pi_i)=\D (\lambda)(\pi_i)=\pi_i \circ \lambda$. That is, $\sigma(\pi_i)=\lambda_i$, which settles the claim in both directions.  For more details, see \cite[Theorem 42]{DMVNotes}.
\end{proof}
 
Consider now the subcategory $\mathbf{DMV_{fp}^{c^\prime}}$  whose 
objects are algebras of the form $DMV_n\mid_P$ for some rational polyhedron $P$, and whose morphisms are maps $\sigma:DMV_n\mid_P\to DMV_m\mid_Q$ such that $\sigma(\{\pi_1, \dots , \pi_n \})\subseteq MV_m$. This category is equivalent to $\mathbf{MV_{fp}}$, since the above condition entails that both categories are dual to $\mathbf{RatPol}_{[0,1]}^{\mathbb Z}$.  Similarly, let $\mathbf{RMV_{fp}^{c^\prime}}$ be a category that has objects of the form $RMV_m\mid_P$, with $P$ a rational polyhedron and arrows  $\sigma:RMV_n\mid_P\to RMV_m\mid_Q$ such that $\sigma(\{\pi_1, \dots , \pi_n \})\subseteq MV_m$. Hence $\mathbf{RMV_{fp}^{c^\prime}}$, the opposite of $\mathbf{RatPol}_{[0,1]}^{\mathbb Z}$ and  $\mathbf{MV_{fp}}$ are equivalent. In this setting the MV-elements are the elements of the MV-algebra generated by the restriction of the projection maps in $DMV_n\mid_P$ or $RMV_n\mid_P$. Once again, our maps are asked to preserves MV-elements.

\subsection{\L -generated logical theories}\label{L-theories}
\noindent The approaches presented in the previous subsections where motivated by the fact that, in order to enlighten subcategories of finitely presented DMV-algebras and Riesz MV-algebras that are equivalent to finitely presented MV-algebras, one has to characterize the \emph{MV-elements} of the algebra at hand. This is a relatively straightforward task when provided with a concrete representation of the categories, but it is difficult  in an abstract setting. The characterization of MV-elements is  important for  defining the right kind of morphisms for the subcategories considered. Indeed, in every approach,  we fundamentally used morphisms that  preserve MV-elements.

On the level of objects, the situation was reflected on ideals:  in order to obtain categorical equivalences for $\mathbf{MV_{fp}}$, the ideal that presents the algebra has to be {\em MV-generated}, which turns into \emph{being generated by an element of the free MV-algebra}. This can be said, of course, once we have fixed a presentation $DMV_n/I$ of the DMV-algebra considered and we have embedded $MV_n$ in $DMV_n$ (the same applies for Riesz MV-algebras). Remark \ref{rem:idRMV} allowed us to prove that any principal ideal in $DMV_n$ is MV-generated (see Lemma \ref{lem:idealsRMV}), while this is no longer true for Riesz MV-algebras, due to Proposition \ref{nonmv}.%, which entails that not every principal ideal of a free Riesz MV-algebra is MV-generated.

This remark, together with the correspondence between theories in logic and filter in the Lindenbaum-Tarski algebra, justifies the notion of {\em \L -generated theory}.  If $\Theta$ be a set of formulas in \L ukasiewicz logic, then the set $Thm(\Theta, \mbox{\L})=\{ \varphi \in Form_{L}\mid \Theta \vdash_{\L} \varphi\}$ the  theory determined by $\Theta$ in \L . Similar definitions are considered for $\mathcal{Q}$\L\ and $\mathcal{R}$\L .

\bdfn If $T$ is a theory in $\mathcal{Q}$\L\ then 
we say that $T$ is {\em \L -generated} if there exists a set of formulas $\Theta\subseteq Form_{L}$  such that $T=Thm(\Theta, \mathcal{Q}\mbox{\L})$. The \L -generated theories of $\mathcal{R}$\L\ are defined similarly.
\edfn

In \cite{DinLeu} it is proved that deductive systems corresponds to filters (and then ideals) in the Lindenbaum-Tarski algebra. We can easily prove the same in the logics $\mathcal{Q}$\L\ and $\mathcal{R}$\L . Moreover, if  $\Theta=\{ \varphi \}$, the ideal is the principal ideal generated by the function that corresponds to $\varphi$. 

We recall that the correspondence is defined as follows:

(i) for any filter $F$ in $MV_n$, $\{ \varphi \mid f_{\varphi}\in F\}$ is a deductive system,

(ii) for any deductive system $\Theta$, the set $\{f_{\varphi}\mid \varphi \in \Theta \}$ is a filter of $MV_n$, for an appropriate $n$.

Thus, the following proposition is straightforward.

\bprop \label{prop:theories}
Let $\varphi$ be a formula of $\mathcal Q$\L . There exists a formula $\beta$ of \L \ such that  $Thm(\varphi, {\mathcal Q}\mbox{\L})= Thm(\beta, {\mathcal Q}\mbox{\L})$.
\eprop
\begin{proof}
We shall assume that all formulas involved have at most $n$ variables. It follows from Remark \ref{lem:idealsRMV} that, if $f_{\varphi}$ is the formula that correspond to $\varphi$, there exists a $g\in MV_n$ such that $(f_{\varphi}]_{DMV}=(g]_{DMV}$. Moreover, since $DMV_n$ is also isomorphic with the algebras of term functions of $\mathcal{Q}$\L , there exists $\beta \in Form_{\mathcal{Q}\L}$ such that $g=f_{\beta}$. Thus, we have the following
\begin{center}
$Thm(\varphi, {\mathcal Q}\mbox{\L})\leftrightarrow (f_{\varphi}]_{DMV}=(f_{\beta}]_{DMV}\leftrightarrow Thm(\beta, {\mathcal Q}\mbox{\L})$,
\end{center}
and the claim is settled.
\end{proof}

As a consequence of the above result, any theory in  $\mathcal{Q}$\L\  is \L -generated. Note that we cannot prove an analogous of Proposition \ref{prop:theories}  for Riesz MV-algebras. A concrete counterexample is given by the formula $\varphi=\neg \left(\Delta_{\frac{\sqrt{2}}{2}}v_1 \oplus \Delta_{\frac{\sqrt{2}}{2}}v_1 \right)$, Indeed, with the notation of Proposition \ref{prop:theories}, $f_{\varphi}$ is the function $f:[0,1]\to [0,1]$ from Example \ref{example:nonmv}. Since $(f]_{RMV}$ cannot be generated by an element of $MV_1$ or $DMV_1$, there isn't a $\beta\in\textbf{\L}$ such that $Thm(\varphi, {\mathcal R}\mbox{\L})= Thm(\beta, {\mathcal Q}\mbox{\L})$.

Whence, in order to obtain something in this vein for Riesz MV-algebras, one should be able to give a precise characterization of those ideals of the form $(f]_{RMV}$, where $f$ is a piecewise linear function with at least one irrational coefficient, in terms of piecewise linear functions with rational coefficients. One possible approach relies on the results on limits of formulas proved in Section \ref{seclogic}.

\bdfn (i) 
An ideal $I$ of $RMV_n$, $n\in \mathbb{N}$, is said to be \emph{norm-closed} if, whenever $f_1, f_2, \dots ,f_m, \dots $ is a sequence of elements of $I$ and $\{ f_m\}_{m\in \mathbb{N}}$ uniformly converges to $f$, then $f\in I$. A similar definition can be given for filters.

(ii) A $\sigma$-ideal $I$ of $RMV_n$ is an ideal closed for all existing countable suprema. That is, if $\{ f_m \}_{m\in \mathbb{N}}\subseteq I$ and $f=\sup_m f_m$ exists in $RMV_n$, then $f\in I$.

(iii) A $\sigma$-filter $F$ of $RMV_n$ is a filter closed for all existing countable infima. That is, if $\{ f_m \}_{m\in \mathbb{N}}\subseteq F$ and $f=\inf_m f_m$ exists in $RMV_n$, then $f\in F$.

We recall that any $\sigma$-ideal is norm-closed, see \cite[Theorem 83.22]{RSZan} and subsequent discussion.
\edfn

The definition of a norm-closed ideal, well known in the theory of Riesz spaces, suggested us the following deduction rule.
$$
\text{($\star$)\quad \quad if } \varphi=\lim_m\varphi_m \text{\quad then \quad} \frac{\varphi_1, \varphi_2, \dots, \varphi_m,\dots}{\varphi}
$$

\bfact
Since the syntactical definition of limit is equivalent with uniform convergence, the limit is continuous with respect to the MV-algebraic operations. Thus, if $\lim_n\varphi_n=\varphi$, then $\lim_n \neg \varphi_n=\neg \varphi$. This entails that we can analogously define \emph{norm-closed filters} and they are again dual to norm-closed ideals. Moreover, it is easily seen that $\sigma$-filters and $\sigma$-ideals are dual to each other.
\efact

\bdfn
$\mathcal{R}$\L$^*$ is the logic obtained from $\mathcal{R}$\L\ by adding the deduction rule ($\star$).
\edfn

\bprop \label{pro:RL*}
The theories of $\mathcal{R}$\L$^*$ are in correspondence with norm-closed ideals of the Lindenbaum-Tarki algebra of $\mathcal{R}$\L .
\eprop

\begin{proof}
We need to prove that, for any $T=Thm(\Theta, \mathcal{R}\text{\L}^*)$ the set $\{ f_{\psi} \mid \psi\in T\} $ is a norm closed filter, and for any norm closed filter $F$ the set $\{ \varphi \mid f_{\varphi}\in F\}$ is a deductive system. The general correspondence between  theories and filters of the Lindenbaum-Tarski algebra is well-known \cite[Proposition 6.3.15]{DinLeu}. We only need to check that ($\star$) implies that the corresponding filter is norm-closed and viceversa. But this is straightforward by definition.
\end{proof}

\bprop
Let $\varphi$ be a formula of $\mathcal{R}$\L . There exists a sequence of formulas $\Theta=\{\varphi_n\}_{n\in \mathbb{N}}\subseteq Form_{\mathcal{Q}L}$ such that $Thm(\varphi, \mathcal{R}\text{\L}^*)= Thm(\Theta, \mathcal{R}\text{\L}^*)$.
\eprop
\begin{proof}

By  Theorem \ref{imp}, there exists a decreasing sequence 
$\{\varphi_n\}_{n\in \mathbb{N}}$ such that $\varphi =\lim_n \varphi_n$. We set 
$\Theta=\{ \varphi_n\mid n\in {\mathbb N} \}$.
One inclusion is direct consequence of  the definition of $\mathcal{R}$\L$^*$.  In order to prove the other inclusion, we note that $f_{\varphi_n}$ converges uniformly to $f_{\varphi}$ and  $f_{\varphi}=\inf_n f_{\varphi_n}$ (cfr \cite[Theorems 16.1 and 16.2]{RS}). Since the set $F=\{ f_{\psi}\mid \psi \in Th(\varphi, \mathcal{R}\text{\L}^*) \}$ is a filter and clearly $f_{\varphi}\in F$, we get that $f_{\varphi_n}\in F $ for any $\varphi_n \in \Theta$, which implies that $\varphi_n$ belongs to $Thm(\varphi, \mathcal{R}\text{\L}^*)$ for any $n$ and the claim is settled.
\end{proof}

We analyzed the conditions that allow us to establish categorical equivalences between subcategories of finitely presented MV-algebras, DMV-algebras and Riesz MV-algebras. 
We have worked with hulls (via tensor product), categories of presentation (using ideas from the theory of institutions) and restriction of free algebras to polyhedra (following the Baker-Beynon duality). All these approaches are connected by one issue: the inability to characterize a peculiar MV-subalgebra of a finitely presentable DMV-algebra or  Riesz MV-algebra. Such elements are called  \emph{MV-elements}, a clear description of them being possible only when we choose a concrete way to define the finitely presentable algebras.  Consequently we have defined the \emph{MV-generated ideals}, which are  ideals in a DMV-algebra or  Riesz MV-algebra that are generated by MV-elements. This leads us to an analysis of those logical theories in the expansions of \L ukasiewicz logic \L \ that are  generated by formulas of \L . Finally,  we characterize the \L -generated theories of ${\mathcal R}$\L \ using the  notion of syntactical limit. 

\subsection*{Acknowledgement}
\noindent %A. Di Nola and S. Lapenta acknowledge partial support from the Italian National Research Project (PRIN2010-11) entitled \textit{Metodi logici per il trattamento dell’informazione}.
I. Leu\c stean was supported by a grant of the Romanian National Authority for Scientific Research and Innovation, CNCS-UEFISCDI, project number PN-II-RU-TE-2014-4-0730. She is also grateful  to  Traian \c Serb\u anu\c t\u a for enlightening conversations on the general theory of presentations and finitely presented structures. 
%The authors are deeply grateful to the two anonymous referees whose comments and suggestions had a major impact on our work.

\end{document}